\documentclass{article}
\usepackage{a4,bbold,pstricks}

\usepackage{amsopn}

\newtheorem{theorem}{Theorem}

\newtheorem{corollary}{Corollary}
\newtheorem{lemma}{Lemma}

\DeclareSymbolFont{Extrasymb}{U}{msa}{m}{n}
\DeclareMathSymbol\square\mathrel{Extrasymb}{"03}
\def\qed{{\leavevmode\unskip\nobreak\hfil\penalty 50\hskip 1em%
  \hbox{}\nobreak\hfil\lower 1pt\hbox{$\square $\kern-.5pt}\parfillskip 0pt
  \finalhyphendemerits 0\par}}

\newcommand{\nc}{\newcommand}
\nc{\ep}{\vspace{0.3cm}}
\nc{\Tr}{\mbox{\bf Tr}}
\nc{\nb}{\nabla}
\nc{\si}{\sigma}
\nc{\Si}{\Sigma}
\nc{\sid}{\si_D}

\nc{\LL}{{\cal L}}
\nc{\HH}{\mathcal{H}}
\nc{\SiM}{\Sigma M}
\nc{\SiN}{\Sigma N}
\nc{\RSi}{R^\Sigma}
\nc{\nbM}{\nabla^M}
\nc{\nbN}{\nabla^N}
\nc{\End}{\mbox{\rm End}}
\nc{\Hom}{\mbox{\rm Hom}}
\nc{\<}{\left\langle}

\nc{\Spin}{\mbox{\rm Spin}}
\nc{\SO}{\mbox{\rm SO}}
\nc{\spin}{\mbox{\goth spin}}
\nc{\so}{\mbox{\goth so}}
\nc{\noi}{\noindent}
\nc{\DM}{D^M}
\nc{\DN}{D^N}
\nc{\DB}{\mbox{\goth D}^B}
\nc{\Mu}{\mbox{\goth M}}
\nc{\Nu}{\mbox{\goth N}}
\nc{\bis}{,\ldots,}
\nc{\muhut}{{\hat{\mu}}}
\nc{\nbmuhut}{\nabla^\muhut}
\nc{\vp}{\varphi}
\nc{\del}{\partial}
\nc{\eh}{\frac{1}{2}}
\nc{\E}{\mathcal E}
\nc{\udot}{\dot{\cup}}
\nc{\bigudot}{\operatorname*{\stackrel{\raisebox{0pt}{.}}{\bigcup}}}
\nc{\uij}{_{i,j}}
\nc{\Id}{\mbox{Id}}
\nc{\lm}{{L^2(M)}}
\nc{\epnoi}{\ep\noi}
\nc{\Cu}{C^\infty}
\nc{\Sjn}{\sum_{j=1}^n}
\nc{\Sj}{\sum_{j}}
\nc{\NN}{\mathcal{N}}
\nc{\punkt}{{\hspace{-0.17cm}\bf .}}{}
\nc{\vol}{\mbox{vol}}
\nc{\FF}{\mathcal{F}}
\nc{\DD}{\mathcal{D}}

\newfont{\Bbb}{msbm10}
\newfont{\goth}{eufm10}

\begin{document}

\title{Zero Sets of Solutions to Semilinear Elliptic Systems of First Order}
\author{Christian B\"ar}
\date{28.~April 1998}
\maketitle

\begin{abstract}
\noindent
Consider a nontrivial solution to a semilinear elliptic system of first
order with smooth coefficients defined over an $n$-dimensional manifold.
Assume the operator has the strong unique continuation property.
We show that the zero set of the solution is contained in a countable
union of smooth $(n-2)$-dimensional submanifolds.
Hence it is countably $(n-2)$-rectifiable and its Hausdorff dimension is 
at most $n-2$.
Moreover, it has locally finite $(n-2)$-dimensional Hausdorff measure.
We show by example that every real number between 0 and $n-2$ actually
occurs as the Hausdorff dimension (for a suitable choice of operator).
We also derive results for scalar elliptic equations of second order.



\end{abstract}


\setcounter{section}{-1}
\section{Introduction}

Many important geometric objects are defined as zero sets of solutions
of certain elliptic differential operators of first order.
The most prominent classical example is provided by algebraic geometry.
Analytic varieties are zero sets of holomorphic functions which in turn
are characterized by the Cauchy-Riemann equations.
More recently zero sets of spinors satisfying a Dirac equation have 
become important.
In \cite{taubes96a} pseudoholomorphic curves in symplectic 4-manifolds
have been constructed using zero sets of harmonic spinors.
Spinors have become a tool to construct conformal immersions of surfaces 
in $\R^3$ with prescribed mean curvature.
Here the zeros are the bad points where the construction does not work.
All this indicates that a general study of the zero set of solutions of 
first order elliptic equations should be useful and important.

Studying the structure of the zero set $\NN$ of such a solution splits
into two quite different problems.
To see this decompose $\NN$ into the set of zeros of finite order, $\NN_{fin}$,
and those of infinite order, $\NN_\infty$.
There is a vast literature concerned with $\NN_\infty$,
the upshot being that it is typically empty.
An operator with $\NN_\infty = \emptyset$ for all its nontrivial solutions is 
said to have the strong unique continuation property.
One knows classical criteria which ensure that an operator has this property, 
see e.g.\ \cite[Ch.~IX]{egorov86a}.
It seems that all operators appearing ``naturally'' in geometry are of this
type \cite{kazdan88a}.
For example Dirac operators (in the most general sense) fall into this class.

Elliptic operators of first order do not always have the
strong unique continuation property, however.
Here is an example.
It is the reduction to first order of the second order
example given in \cite[Ex. 1.11]{kazdan88a}.

\epnoi
{\bf Example.}
Let $M=\R^2$.
Let $\Delta_1 = -\frac{\del^2}{\del x^2} - \frac{\del^2}{\del y^2}$ be the 
Laplace operator and let
$\Delta_2 = -\frac{\del^2}{\del x^2} - \frac{1}{2}\frac{\del^2}{\del y^2}$,
both operators acting on $\Cu(\R^2,\R)$.
The fourth order operator $\Delta_1 \Delta_2$ satisfies the
assumptions of a theorem by Alinhac \cite{alinhac80a} which
tells us that there exist functions $u,a \in \Cu(U,\R)$, defined
on a neighborhood $U \subset \R^2$ of $0$, vanishing
of infinite order at $0$, but not identically zero, and solving
$$
\Delta_1 \Delta_2 u = a \cdot u .
$$
Now let $D_1 = d + \delta$ act on $\Cu(\R^2,\bigoplus_{j=0}^2 
\Lambda^j T^\ast\R^2)$.
This operator is a square root of $\Delta_1$ in the sense that 
$D_1^2|_{\Cu(\R^2,\Lambda^0 T^\ast\R^2)} = \Delta_1$.
Similarly, let $D_2$ acting on $\Cu(\R^2,\bigoplus_{j=0}^2 
\Lambda^j T^\ast\R^2)$ be a square root of $\Delta_2$.
Then we have 
$$
D_1D_1D_2D_2 u = au .
$$
Therefore $\vp = (u,D_2u,D_2^2u,D_1D_2^2u)$ solves the linear elliptic
system
\begin{equation}
D\vp =
\left(
\begin{array}{cccc}
D_2 & -1 & 0 & 0 \\
0 & D_2 & -1 & 0 \\
0 & 0 & D_1 & -1 \\
-a & 0 & 0 & D_1
\end{array}
\right)
\vp = 0
\label{nonstrongbsp}
\end{equation}
Since $u$ has a zero of infinite order at 0 the same is true for $\vp$.
Therefore the linear elliptic operator $D$ does not have the strong 
unique continuation property.
Since every ``component'' in the matrix in (\ref{nonstrongbsp}) acts 
itself on sections of a 4-dimensional vector bundle, (\ref{nonstrongbsp})
amounts to a system of 16 equations when spelled out.

\epnoi
Such counterexamples always seem somewhat artifial.
In interesting cases the strong unique continuation property is usually
satisfied.
Therefore we focus our attention on the other component of the zero set, 
$\NN_{fin}$.
Unless the operator has analytic coefficients or the underlying
manifold is of dimension $n \le 2$ surprisingly little seems to be known 
about it.
Classical theorems on uniqueness in the Cauchy problem tell us that 
$\NN_{fin}$ does not contain certain hypersurfaces.
But this does not really mean much, $\NN_{fin}$ could still be a very
irregular set of any Hausdorff dimension.

It turns out that the zero sets are in general very irregular but also
very well-behaved depending on the point of view.
We will show by example in Theorem \ref{wild} that every closed
subset of a submanifold of codimension 2 is the zero set of a solution
for some first order elliptic operator.
This operator has the strong unique continuation property.
Hence the Hausdorff dimension of the zero set can be any {\em real}
number between 0 and $n-2$ where $n$ is the dimension of the underlying
manifold.

On the other hand, the zero sets have the following regularity properties.
Our Main Theorem says that they are contained in a countable union
of smooth submanifolds of codimension 2.
In particular, they are countably $(n-2)$-rectifiable and the
Hausdorff dimension is at most $n-2$.
Moreover, they have locally finite $(n-2)$-dimensional Hausdorff measure.
If the underlying manifold is a surface, $n=2$, then $\NN_{fin}$ is discrete.
We give an upper bound for the Hausdorff measure in small balls in terms
of the vanishing order of the solution.
All this shows that one can apply methods from geometric measure theory
to such zero sets.
We hope that this will be useful in the future.

Our main tool is the Malgrange Preparation Theorem which states that 
the zero set of a smooth real-valued function vanishing of finite order
can locally be described as the zero set of another function which is 
polynomial in one of its variables.
>From this one gets countable $(n-1)$-rectifiability fairly easily.
The difficulty in the proof is to show that the zero sets of the 
various components of our solution intersect in such a way that they
decrease the dimension once more.
This requires an algebraic study of certain ``resultants''.
This approach has first been used by the author in \cite{baer97b} to
investigate zero sets of Dirac operators.
In this paper we extend those results in two respects.
Firstly, we enlarge the class of operators from (linear) Dirac
operators to semilinear elliptic operators of first order.
This seems to be the largest natural class of operators for which our results
can be expected to hold.\footnote{Actually, this is not quite true.
Once our Main Theorem is shown to hold for semilinear equations it is
clear it also holds for the larger class of quasilinear equtions.
The obvious details are left to the reader.}{}
Secondly, we get more precise information about the zero sets such as
bounds on their Hausdorff measure.

By a suitable reduction of order we obtain corollaries for scalar elliptic
equations of second order as well.
Decompose the zero set $\NN$ of a solution, often also called its nodal 
set, into a smooth hypersurface and into the critical zero set $\NN_{crit}$ 
where the function and its gradient vanish simultaneously.
We give simple new proofs of the following facts:
$\NN$ is countably $(n-1)$-rectifiable with locally finite $(n-1)$-dimensional
Hausdorff measure and $\NN_{crit}$ is countably $(n-2)$-rectifiable with 
locally finite $(n-2)$-dimensional Hausdorff measure.
Again, we get upper bounds for the Hausdorff measure in small balls.
These facts have been obtained by various authors over the time, some
of them only very recently.
It is amusing that they follow easily from our results on first order systems
by a simple reduction of order.

The paper is organized as follows.
In the first section we give the necessary definitions and a few
examples.
We state the Main Theorem and deduce some immediate corollaries.
The second section contains the proof of the Main Theorem.
In the third section we construct examples with very irregular zero sets.
The last section contains the discussion of scalar equations of second order.

All differential operators in this paper are assumed to have $\Cu$-coefficients.

\epnoi
{\bf Acknowledgements.}
It is a pleasure to thank J.~Dodziuk, H.~Kalf, P.~Lax, and F.H.~Lin for 
valuable discussion and helpful hints.
This paper was written while the author enjoyed the hospitality of the
Courant Institute of Mathematical Sciences.


\section{Semilinear Elliptic Differential Operators of First Order}

Let $M$ be an $n$-dimensional connected differential manifold, 
let $E,F \to M$ be real vector bundles over $M$.
The case of complex vector bundles is included in our discussion
because we can simply forget the complex structures and regard
the bundles as real vector bundles.
Let $D : \Cu(M,E) \to \Cu(M,F)$ be a linear differential operator of first
order. 
Here $\Cu(M,E)$ and $\Cu(M,F)$ denote the spaces of smooth sections
in the bundles $E$ and $F$ resp.
If we introduce local coordinates $x_1\bis x_n$ on $M$ and trivialize
the bundles, then $D$ takes the form
$$
D = \sum_{j=1}^n A_j(x)\frac{\del}{\del x_j} + B(x)
$$
where $A_j$ and $B$ are smooth matrix-valued functions.

The {\em principal symbol} $\si_D$ of $D$ associates to each covector
$\xi \in T^\ast_x M$ with base point $x$ a homomorphism $\si_D(\xi) : 
E_x \to F_x$ which is characterized as follows:

Choose a smooth function $f$ defined in a neighborhood of $x$ such
that $f(x)=0$ and $df(x)=\xi$.
Take an arbitrary $\vp \in E_x$ and extend it smoothly to a section
$\Phi$ of $E$ in a neighborhood of $x$.
Then
$$
\si_D(\xi)\vp = D(f\Phi)(x).
$$
It is easy to see that $\si_D(\xi)$ does not depend on the choices of
$f$ and $\Phi$.
In coordinates, for $\xi = \sum_{j=1}^n \xi_j dx_j$, this homomorphism
is given by the matrix 
$$
\si_D(\xi) = \Sjn A_j(x)\xi_j.
$$
Very frequently in the literature there is an additional factor of $i$
in the definition of the principal symbol.
But since we are dealing with real bundles and operators introducing
a factor of $i$ would seem somewhat artificial.
The only disadvantage of our convention is a minus sign in the formula
for the symbol of the formally adjoint operator, $\si_{D^\ast}(\xi) =
-\si_D(\xi)^\ast$.
But this will be of no relevance to our discussion.

Note that the map $T^\ast_xM \to \Hom(E_x,F_x)$, $\xi\mapsto\si_D(\xi)$, 
is linear.
This is a special property of first order operators.

For a section $\vp\in\Cu(M,E)$ and a function $f\in\Cu(M,\R)$ we have the
formula
$$
D(f\vp) = fD(\vp) + \si_D(df)\vp .
$$

The operator $D$ is called {\em elliptic} if the symbol $\si_D(\xi)$
is an isomorphism for all $\xi\not= 0$.
In particular, $E$ and $F$ must then have the same rank.

Now let $V:E \to F$ be a smooth fiber-preserving map, i.e.\
$V(E_x) \subset F_x$ for all $x \in M$.
We say that $V$ {\em respects the zero section} if $V(0)=0$ for $0\in E_x$
and all $x$.
Of course, this is automatic if $V$ is fiberwise linear.

An operator of the form
$$
L=D+V : \Cu(M,E) \to \Cu(M,F)
$$
where $D$ is a linear first order differential operator is called
a {\em semilinear} differential operator of first order.
We define $\si_L := \si_D$ as the principal symbol of $L$.
We say $L$ is {\em elliptic} if $D$ is and we say $L$ {\em respects the
zero section} if $V$ does.

Note that the decomposition $L=D+V$ of a semilinear operator is unique
up to linear zero-order terms.
Hence all concepts are well-defined.

For a section $\vp\in\Cu(M,E)$ we define its {\em zero set} 
$$
\NN(\vp) = \{ x\in M\ |\ \vp(x)=0 \}.
$$
We say that $\vp$ {\em vanishes at $x$ of order $k$} if when expressed in
a local trivialization all derivatives up to order $k-1$ of all
components $\vp_i$, $\vp = (\vp_1\bis \vp_N)$, vanish at $x$,
$$
\frac{\del^m\vp_i}{\del x_{j_1}\cdots\del x_{j_m}}(x) = 0,
$$
for $1\le i \le N$, $1 \le j_\nu \le n$, and $0 \le m \le k-1$.
This condition is independent of the choice of coordinates and
local trivialization.

We split the zero set into two parts,
$$
\NN(\vp) = \NN_{fin}(\vp) \udot \NN_\infty(\vp),
$$
where $\NN_{fin}(\vp)$ is the set of zeros of $\vp$ of finite
order, i.e.\ the set of those zeros for which there exists $k \in \N$ 
such that $\vp$ vanishes of order $k$ but not of order $k+1$.
Accordingly, $\NN_\infty(\vp)$ is set of zeros of infinite order.

A differential operator $L$ respecting the zero section is said 
to have the {\em strong unique continuation property} if for all local 
solutions $\vp$ of $L\vp=0$
either $\NN_\infty(\vp) = \emptyset$ or $\vp\equiv 0$.

If $L$ is a semilinear elliptic differential operator of first
order on a connected 1-dimensional manifold, then the uniqueness theorem
for ordinary differential equations tells us that $\NN(\vp) =
\emptyset$ unless $\vp \equiv 0$.
In particular, $L$ has the strong unique continuation property.
In dimension $n\ge 2$, the strong unique continuation
property sometimes fails even for linear elliptic operators of first
order as we have seen in the introduction.
There is a vast literature on conditions on differential operators
which imply the strong unique continuation property.
Fortunately, operators arising from geometric problems always seem 
to satisfy it \cite{kazdan88a}.

\epnoi
{\bf Example.}
Let the underlying manifold $M$ carry a Riemannian metric $g$. 
If the principal symbol $\si_D$ of a linear differential operator 
$D:\Cu(M,E)\to\Cu(M,E)$ of first order satisfies the {\em Clifford relations}
$$
\si_D(\xi)\circ\si_D(\eta) + \si_D(\eta)\circ\si_D(\xi) 
+ 2g(\xi,\eta)\cdot\Id_E =0
$$
for all $\xi,\eta \in T^\ast_x M, x \in M$, then $D$ is called a 
{\em Dirac operator}.
If, in addition, $V$ respects the zero section, then we call 
$L = D+V$ a {\em semilinear Dirac operator}.

Such a Dirac operator is certainly elliptic because the Clifford 
relations imply $\si_D(\xi)^2 = -|\xi|^2\cdot\Id_E$.

Moreover, 
solutions $\vp$ of a semilinear Dirac equation $L\vp=0$ locally satisfy 
a differential inequality 
$$
|D^2 \vp| = |D(V(\vp))| \le   C\cdot (|\vp| + |\nb\vp|).
$$
Since $D^2$ is an elliptic differential operator of second order
with scalar symbol Aronszajn's theorem \cite{aronszajn57a} applies 
and tells us that $L$ has the strong unique continuation property.

\epnoi
{\bf Example.}
Let $M$ be an oriented surface equipped with a spin structure.
Let $D:\Cu(M,\SiM)\to\Cu(M,\SiM)$ be the Dirac operator acting on spinors 
and let $H:M\to\R$ be a smooth function.
Define $V:\SiM\to\SiM$ by 
$$
V(\vp)=-H |\vp |^2\vp.
$$
Then $L=D+V$ is a semilinear Dirac operator (which respects the zero section).

This operator has attracted much attention in recent years because
the solutions of $L\vp=0$ give rise to conformal immersions of the 
universal cover of $M-\NN(\vp)$ into $\R^3$ with mean curvature $H$.
This can be regarded as a generalization of the classical Weierstrass
representation of minimal surfaces, see e.g.\ \cite{friedrich97ppa, 
kamberov-pedit-pinkall96ppa,kusner-schmitt96ppa,taimanov97ppa}.

The operator $L$ has the strong unique continuation property. 
Corollary \ref{dirac} to our Main Theorem then says that $\NN(\vp)$ is 
discrete.
Thus $\vp$ defines a multivalued immersion of $M$ into $\R^3$
branched along the discrete set $\NN(\vp)$.

\ep
In this paper we study the set $\NN_{fin}(\vp)$ for solutions
of arbitrary semilinear elliptic operators of first order which
respect the zero section.
For operators having the strong unique continuation property
such as semilinear Dirac operators this means that we will be able to
to control all of $\NN(\vp)$.

\ep
Recall that a subset of $\R^n$ is called {\em countably $k$-rectifiable}
if it can be written as a countable union of images under
Lipschitz maps of bounded closed subsets of $\R^k$, c.f.~\cite{federer69a}.
A subset of a manifold is countably $k$-rectifiable if it is so
in coordinate charts.

In fact, we will prove something stronger rather than countable
rectifiability.
Therefore we make the following 

\epnoi
{\bf Definition.}
A subset of a differential manifold is called {\em countably 
$k$-$\Cu$-rectifiable} if it is contained in a countable union
of smooth $k$-dimensional submanifolds.

\epnoi
Of course, a set which is countably $k$-$\Cu$-rectifiable is also
countably $k$-recti\-fiable.

\epnoi
Also recall the definition of Hausdorff measure density 
\cite[2.10.19]{federer69a}.
Let $\alpha(m)$ denote the $m$-dimensional volume of the 
unit ball in $\R^m$.
Denote the $m$-dimensional Hausdorff measure by $\HH^m$.
If $N$ is a subset of a Riemannian manifold and $p \in N$, 
then the limit 
\begin{equation}
\Theta^{\ast m}(N,p) = \limsup_{r \searrow 0} 
\frac{\HH^m(N \cap B(p,r))}{\alpha(m)r^m} \in [ 0,\infty ]
\label{density}
\end{equation}
is called {\em $m$-dimensional upper Hausdorff density} of $N$ at $p$.

If for example $N$ is an $m$-dimensional submanifold, then
$\Theta^{\ast m}(N,p) = 1$.

\epnoi
{\bf Main Theorem.}
{\em Let $M$ be a connected $n$-dimensional differential manifold.
Let $L$ be a semilinear elliptic differential operator defined over $M$
which respects the zero section.
Let $\vp\not\equiv 0$ satisfy $L\vp=0$.

Then $\NN_{fin}(\vp)$ is a countably $(n-2)$-$\Cu$-rectifiable set.
At each point $p \in \NN_{fin}(\vp)$ the $(n-2)$-dimensional upper 
Hausdorff density has a bound
$$
\Theta^{\ast n-2}(\NN(\vp),p) \le C(n)k^3
$$
where $C(n)$ is a constant depending only on the dimension $n$
and $k$ is the order of vanishing of $\vp$ at $p$.

In particular, we have for the Hausdorff dimension
$$
\dim(\NN_{fin}(\vp)) \le n-2
$$
and if $n=2$, then $\NN_{fin}(\vp)$ is a discrete set.}

\epnoi
The proof will deliver an explicit value for $C(n)$.
For example, we can take $C(n) = \frac{2^{n-3}n(n-1)}{\alpha(n-2)}$.
This constant is probably not optimal.

The condition that $L$ respect the zero section is obviously
necessary for the theorem to hold as one can see already in the 
1-dimensional case.
There are many solutions of ordinary differential equations of first
order not respecting the zero section which do have zeros.
Of course, if the ordinary differential equation respects the zero
section, i.e.\ it is of the form $\frac{d}{dt}u(t) + F(t,u(t)) = 0$
with $F(t,0)=0$, then since $u_0 \equiv 0$ is a solution, a
nontrivial solution $u$ cannot have any zeros by the uniqueness
theorem.

In \cite{baer97b} it was shown that the bound $n-2$ in the Main Theorem
is sharp already in the class of linear Dirac operators.
In Section 3 we will see by example that $\NN_{fin}(\vp)$ can be 
very irregular and can have any {\em real} number $d\in [0,n-2]$
as its Hausdorff dimension.

Combining the Main Theorem with the strong unique continuation
property yields the following corollary which applies in particular
to semilinear Dirac operators.

\begin{corollary}\punkt
\label{dirac}
Let $M$ be a connected $n$-dimensional differential manifold.
Let $L$ be a semilinear elliptic differential operator defined over $M$
which respects the zero section and which has the strong unique continuation
property.
Let $\vp\not\equiv 0$ satisfy $L\vp=0$.

Then $\NN(\vp)$ is countably $(n-2)$-$\Cu$-rectifiable and has locally finite
$(n-2)$-dimensional Hausdorff measure.
In particular, we have for the Hausdorff dimension
$$
\dim(\NN(\vp)) \le n-2
$$
and if $n=2$, then $\NN(\vp)$ is a discrete set.
\qed
\end{corollary}

In the linear case the difference of two solutions is again a solution,
so that the Main Theorem can be regarded as a strong version of the 
uniqueness part in the Cauchy problem.

\begin{corollary}\punkt
\label{unique}
Let $M$ be a connected $n$-dimensional differential manifold.
Let $A \subset M$ be a closed subset of Hausdorff dimension
$\dim A > n-2$.
Let $D$ be a linear elliptic differential operator of first order
defined over $M$ which has the strong unique continuation property.
Let $\vp_1$ and $\vp_2$ satisfy $L\vp_i=0$.

If $\vp_1|_A=\vp_2|_A$, then $\vp_1=\vp_2$ on all of $M$.
\end{corollary}

\noi
Here $A$ replaces the hypersurface in the classical Cauchy problem.

Applications of the Main Theorem to scalar elliptic equations of
second order will be given in Section 4.
We conclude this section with a corollary for a very special but important
elliptic second order system.

\begin{corollary}\punkt
\label{l2form}
Let $M$ be a complete connected $n$-dimensional Riemannian manifold.
Let $\Delta = d\delta + \delta d$ be the Laplace-Beltrami operator
acting on $p$-forms.
Let $\omega \not= 0$ be a square-integrable harmonic $p$-form, 
$\Delta\omega = 0$.

Then $\NN(\omega)$ is countably $(n-2)$-$\Cu$-rectifiable and has locally
finite $(n-2)$-dimensional Hausdorff measure.
In particular, the Hausdorff dimension of $\NN(\omega)$ is at most
$n-2$ and if $n=2$, then $\NN(\omega)$ is discrete.
\end{corollary}

Note that Corollary \ref{l2form} fails if one drops the assumption of 
square-integrability. 
For example, $\omega = x_1 dx_1 \wedge \ldots \wedge dx_p$ is a 
harmonic $p$-form on $\R^n$ but its zero set has codimension 1.
The corollary also fails if one replaces harmonic forms by eigenforms
for positive eigenvalues.

\epnoi
{\em Proof of Corollary \ref{l2form}.}
By \cite[Thm.~26]{rham84a} harmonic $L^2$-forms are closed and coclosed,
$$
(d+\delta)\omega = 0.
$$
The operator $d+\delta$ is a Dirac operator.
Thus it has the strong unique continuation property, $\NN_\infty(\omega)
=\emptyset$, and the Main Theorem gives the result for $\NN(\omega) =
\NN_{fin}(\omega)$.
\qed


\section{The Proof}
In this section we prove the Main Theorem.
We first consider a very special class of differential operators.
Let $X$, $E$, and $F$ be real vector spaces with $\dim(X)=n$ and 
$\dim(E)=\dim(F)=N$.
Let $\si : X^\ast \to \Hom(E,F)$ be a linear map.
This induces a differential operator of first order $\si(\del):
\Cu(X,E) \to \Cu(X,F)$ as follows:
Choose a basis $e_1\bis e_n$ of $X$, let $e_1^\ast\bis e_n^\ast$ 
be the dual basis and put
$$
\si(\del) = \Sjn \si(e_j^\ast)\del_{e_j}.
$$
This definition is easily seen to be independent of the choice
of basis $e_1\bis e_n$.
Such an operator is called an operator {\em with constant coefficients}.
Of course, the principal symbol of $\si(\del)$ is precisely given by 
$\si$.
We now look at the zero set of polynomial solutions of elliptic operators
with constant coefficients.
For elements $(x_1,x_2\bis x_n) \in \R^n$ we use the notation 
$(x_2\bis x_n) = x'$ and denote the projection $x \mapsto x'$ by
$\pi : \R^n \to \R^{n-1}$.

\begin{lemma}\punkt
\label{konstant}
Let $\si(\del): \Cu(\R^n,\R^N) \to \Cu(\R^n,\R^N)$ be an elliptic first
order differential operator with constant coefficients.
Let $\vp = (\vp_1\bis\vp_N) \in \Cu(\R^n,\R^N)$ be a solution of 
$\si(\del)\vp = 0$ such that all 
components $\vp_\nu$ are homogeneous polynomials of the same degree $k$ 
of the form
$$
\vp_\nu(x_1,x') = \alpha_\nu
\left( x_1^k + \sum_{j=0}^{k-1} u_{\nu,j}(x')x_1^j \right) ,
$$
where $\alpha_1\bis\alpha_N \in \R$ and $u_{\nu,j}$ are homogeneous
polynomials of degree $k-j$.
Then:

If $\pi(\NN(\vp)) = \R^{n-1}$, then $\vp \equiv 0$.
\end{lemma}

\noi
{\em Proof.}
We write $\vp$ as a polynomial in $x_1$ with vector-valued coefficients
\begin{equation}
\vp(x) = \sum_{j=0}^k Y_j(x')\cdot x_1^j
\label{phormel1}
\end{equation}
where $Y_k(x') = (\alpha_1\bis\alpha_N)$ and $Y_j(x') = (\alpha_1
u_{1,j}(x') \bis \alpha_N u_{N,j}(x'))$, $j=0\bis k-1$.

We define an elliptic differential operator with constant coefficients by 
$$
D : \Cu(\R^{n-1},\R^N) \to \Cu(\R^{n-1},\R^N),
$$
$$
D = -\si(e_1^\ast)^{-1}\sum_{j=2}^n \si(e_j^\ast)\frac{\del}{\del x_j}.
$$

We compute
\begin{eqnarray*}
0 &=& \si(e_1^\ast)^{-1}\si(\del)\vp \\
&=&
\left( \frac{\del}{\del x_1} - D \right)
\left( \sum_{j=0}^k Y_j(x')\cdot x_1^j   \right) \\
&=&
\sum_{j=0}^{k-1} \left( (j+1)Y_{j+1}(x') - DY_j(x') \right) x_1^j
\end{eqnarray*}

We conclude $Y_{j+1} = \frac{1}{j+1}DY_j$ and hence 
$Y_j = \frac{1}{j!}D^jY_0$.
Thus we can rewrite (\ref{phormel1}) as
\begin{equation}
\vp(x) = \sum_{j=0}^k \frac{1}{j!}D^jY_0(x')\cdot x_1^j
\label{phormel2}
\end{equation}

The assumption $\pi(\NN(\vp)) = \R^{n-1}$ means that for any 
$x' \in \R^{n-1}$ we can find an $x_1(x')\in\R$ such that
\begin{equation}
\sum_{j=0}^k \frac{x_1(x')^j}{j!}D^jY_0(x') = 0
\label{phormel3}
\end{equation}
Roots of polynomials do not depend smoothly on the coefficients of
the polynomial everywhere but smoothness fails only when multiple zeros 
branch to distinct zeros.
Hence we can assume that $x_1(x')$ depends smoothly on $x'$ on a
nonempty open subset of $\R^{n-1}$.
On this subset we apply $D$ to equation (\ref{phormel3}) and obtain
\begin{eqnarray}
0 &=& D\sum_{j=0}^k \frac{x_1(x')^j}{j!}D^jY_0(x') \nonumber \\
&=& \sum_{j=0}^{k-1} \frac{x_1(x')^j}{j!}D^{j+1}Y_0(x') -
    \sum_{j=1}^k \si(e_1^\ast)^{-1}\si\left( d\left(\frac{x_1(x')^j}{j!}
    \right)\right) D^jY_0(x') \nonumber \\
&=& (1-\si(e_1^\ast)^{-1}\si(dx_1(x'))) \sum_{j=0}^{k-1} 
    \frac{x_1(x')^j}{j!}D^{j+1}Y_0(x')
\label{reduziert}
\end{eqnarray}
Now observe
\begin{eqnarray*}
\det[1-\si(e_1^\ast)^{-1}\si(dx_1(x'))] &=&
\det[\si(e_1^\ast)]^{-1} \det[\si(e_1^\ast) - \si(dx_1(x'))] \\
&=&
\det[\si(e_1^\ast)]^{-1} \det[\si(e_1^\ast - dx_1(x'))] \\
&\not=& 0
\end{eqnarray*}
since $e_1^\ast - dx_1(x') \not= 0$ and $\si(\del)$ is elliptic.
Hence (\ref{reduziert}) yields
\begin{equation}
\sum_{j=0}^{k-1} \frac{x_1(x')^j}{j!}D^{j+1}Y_0(x') = 0
\label{phormel4}
\end{equation}
Equation (\ref{phormel4}) is the same as (\ref{phormel3}) with $k$
replaced by $k-1$ and $Y_0$ replaced by $DY_0$.
Repeating this inductively we eventually get
$$
0=D^kY_0(x')=k! Y_k(x') = k! (\alpha_1\bis\alpha_N).
$$
Hence $\vp\equiv 0$.
\qed

\epnoi
To study local properties of differential operators it can be useful 
to approximate the operator in a neighborhood of a given point by an
operator with constant coefficients.
This {\em freezing of coefficients} is done as follows.

Let $L:\Cu(M,E)\to\Cu(M,F)$ be a semilinear differential operator 
which respects the zero section.
Let $p\in M$ be a point.
Then the principal symbol $\si_{L,p} : T^\ast_pM \to \Hom(E_p,F_p)$ 
gives rise to a linear differential operator with constant coefficients
$$
\hat{L}_p = \si_{L,p}(\del) : \Cu(T_pM,E_p) \to \Cu(T_pM,F_p)
$$
defined on the tangent space of $M$ at $p$.

To a section $\vp\in\Cu(M,E)$ and a point $p\in M$ we associate a
homogeneous polynomial $\hat{\vp}_p \in \Cu(T_pM,E_p)$ as follows.
Choose local coordinates near $p$ and look at the Taylor series
expansion of $\vp$ at $p$.
Then $\hat{\vp}_p$ is the homogeneous part of lowest degree in the 
expansion which does not identically vanish,
$$
\vp(p+\xi) = \hat{\vp}_p(\xi) + \mbox{ terms of higher order in }\xi
$$
One checks that this gives rise to a well-defined homogeneous polynomial
on $T_pM$ with coefficients in $E_p$.
E.g., if $\vp(p)\not= 0$, then $\hat{\vp}_p$ is constant, $\hat{\vp}_p(\xi)
=\vp(p)$.
If $\vp$ vansihes of first order at $p$, then $\hat{\vp}_p(\xi)=
\nb_\xi\vp$ for any connection $\nb$ on $E$.
If $\vp$ vanishes of infinite order at $p$, then we set $\hat{\vp}_p
\equiv 0$.

The next lemma tells us how to pass
from an arbitrary solution of a semilinear first order equation to a
polynomial solution of an operator with constant coefficients.

\begin{lemma}\punkt
\label{hauptteil}
Let $M$ be an $n$-dimensional differential manifold.
Let $E$ and $F$ be real vector bundles over $M$.
Let $p\in M$ be a point.
Let $L:\Cu(M,E)\to\Cu(M,F)$ be a semilinear elliptic differential 
operator which respects the zero section.
Then:

If $\vp\in\Cu(M,E)$ satisfies 
$$
L\vp=0,
$$ 
then $\hat{\vp}_p\in\Cu(T_pM,E_p)$ satisfies 
$$
\hat{L}_p \hat{\vp}_p =0.
$$
\end{lemma}

\noi
{\em Proof.}
Expand $\vp$ and the coefficients of $L$ in their Taylor series at $p$,
use $V(\vp) = \mbox{O}(|\vp|)$ and take the lowest order term of $L\vp$.
\qed

\begin{lemma}\punkt
\label{hyperfl}
Let $U \subset \R^n$ be an open neighborhood of $0$.
Let $f : U \to \R$ be a smooth function vanishing of finite order $k$ at $0$,
but not of order $k+1$.
Write $\NN(f) = f^{-1}(0)$ for the zero set.

Then for sufficiently small $r>0$ the set $\NN(f) \cap B(0,r)$
is countably $(n-1)$-$\Cu$-rectifiable and
$$
\Theta^{\ast n-1}(\NN(f),0) \le C(n)k .
$$
\end{lemma}

\noi
Here $C(n)$ is a constant depending only on $n$.
For example, $C(n)=\frac{n2^{n-1}}{\alpha(n-1)}$ works where $\alpha(m)$ is 
the volume of the unit ball in $\R^m$.
See (\ref{density}) for the definition of $\Theta^{\ast n-1}$.

The proof will show a somewhat stronger statement then just countably 
$(n-1)$-$\Cu$-rectifiability of $\NN(f) \cap B(0,r)$.
Namely, for a generic choice of cartesian coordinate system on $\R^n$
the following is true:

Pick any index $j \in \{ 1 \bis n \}$.
Then $B(0,r)$ can be written as a countable union of subsets in which 
the $x_j$-component of the points in $\NN(f)$ is a
smooth function defined on bounded subsets of the hyperplane
perpendicular to the $x_j$-axis.
In slightly other words, near 0, $\NN(f)$ is contained in a countable union 
of smooth graphs over the hyperplane $e_j^\perp$.


\begin{center}
\pspicture(1,0)(14,6)

\rput(11,5.5){Case $n=2$}

\psline{->}(7,0)(7,6)
\rput(7.5,5.5){$x_2$}
\psline{->}(4,3)(10,3)
\rput(9.5,2.5){$x_1$}

\pscurve(5,1.5)(6,2)(7,3)(8,4)(9,4.5)
\pscurve(5.5,1)(6,2.5)(7,3)(8.5,4)(9,5)
\pscurve(5,4.5)(6,4)(7,3)(8,2)(9,1.5)

\rput(6,4){\psframebox*[framearc=0.5]{$\supset\mathcal{N}(f)$}}

\psline{->}(7.2,1.9)(7.9,1.9)
\psline{->}(8.2,2.8)(8.2,2.1)

\endpspicture
\end{center}
\begin{center}
\bf Fig.~1
\end{center}

\epnoi
{\em Proof.}
Since $f$ vanishes of finite order $k$ at $0$, we can choose an orthonormal
basis $e_1 \bis e_n$ of $\R^n$ such that $D^kf(0)(e_i \bis e_i) \not= 0$
for all $i = 1 \bis n$.
We use cartesian coordinates on $\R^n$ with respect to which $e_1 \bis e_n$
are the standard basis.

By Malgrange's Preparation Theorem \cite[Ch.~V]{malgrange66a} we can write
$f$ in a neighborhood of $0$ as 
\begin{equation}
\label{malgrange1}
f(x) = v(x) \left( x_1^k + \sum_{j=0}^{k-1} u_j(x') x_1^j \right)
\end{equation}
where $v, u_0 \bis u_{k-1}$ are smooth functions defined in a neighborhood
of $0$, $v$ nowhere vanishing and $u_j$ vanishing at 0 of order $k-j$ at least.
So $\NN(f)$ is given near $0$ as the zero locus of a polynomial in the
$x_1$-variable of degree $k$.
Hence if we define the $m$-dimensional cube $W^m(0,r) = \stackrel{
m \mbox{ factors}}{\overbrace{[-r,r] \times \cdots \times [-r,r]}}$, then
under the projection onto the $x'$-hyperplane each $x' \in W^{n-1}(0,r)$
has at most $k$ preimages contained in $\NN(f) \cap W^n(0,r)$.

Now recall the following known facts about how roots of polynomials
depend on the coefficients of the polynomial \cite[L.~6]{yomdin84a}:
$\R^k$ can be covered by countably many bounded closed sets
$\R^k = \bigcup_\mu A_\mu$ such that number of real roots
of the polynomial $F_u(t) = t^k + \sum_{j=0}^{k-1} u_j t^j$
is the same number $k_\mu \in \{ 0 \bis k \}$ for all $u = (u_0 \bis 
u_{k-1}) \in A_\mu$.
Moreover, there exist smooth functions $\xi_{\mu,i} : \R^k \to \R$,
$i = 1 \bis k_\mu$, such that $\xi_{\mu,1}(u) < \cdots < \xi_{\mu,k_\mu}(u)$
are precisely the roots of $F_u$ whenever $u \in A_\mu$.

This together with (\ref{malgrange1}) shows that $\NN(f) \cap W^n(0,r)$ is 
contained in the union of the graphs
of the smooth functions $x_1 = \xi_{\mu,i}(u(x'))$, $x' \in W^{n-1}(0,r)$.
Hence $\NN(f)$ is countably $(n-1)$-$\Cu$-rectifiable near $0$.

To estimate the Hausdorff density we note that by the choice of the
coordinate system we can apply Malgrange's preparation theorem to any
of the variables $x_j$, not just $x_1$.
Therefore $\NN(f) \cap W^n(0,r)$ contains for each $x' \in W^{n-1}(0,r)$
at most $k$ preimages under each projection to the hyperplane $\R^{n-1} =
e_j^{\perp} \subset \R^n$.
Now the second estimate in \cite[3.2.27]{federer69a} yields for small $r>0$
$$
\HH^{n-1}(\NN(f) \cap W^n(0,r)) \le nk\cdot\vol(W^{n-1}(0,r)) = 
nk(2r)^{n-1} .
$$
Thus
\begin{eqnarray*}
\Theta^{\ast n-1}(\NN(f),0) 
&=&
\limsup_{r \searrow 0} \frac{\HH^{n-1}(N \cap B(p,r))}{\alpha(n-1)r^{n-1}} \\
&\le&
\limsup_{r \searrow 0} \frac{nk(2r)^{n-1}}{\alpha(n-1)r^{n-1}} \\
&=&
\frac{n2^{n-1}}{\alpha(n-1)} k
\end{eqnarray*}
\qed

\epnoi
With these preparations we can now proceed to the proof of the Main Theorem.

\epnoi
{\em Proof of Main Theorem.}
Let $L:\Cu(M,E) \to \Cu(M,F)$ be a semilinear elliptic operator respecting
the zero section and let $\vp \in \Cu(M,E)$ solve $L\vp=0$.
Let $p\in\NN_{fin}(\vp)$.
We will show that $\NN(\vp)$ is countably $(n-2)$-$\Cu$-rectifiable near $p$.

Since $\vp$ vanishes at $p$ of some order $k$ but not of order $k+1$
we can choose the coordinates near $p=0$ in such a way that the homogeneous 
part of order $k$ has the form
$$
\hat{\vp}_p(x) = a \cdot x_1^k + \mbox{ terms of lower degree in }x_1
$$
with $a\in E_p$, $a\not= 0$.
We trivialize $E$ near $p$ in such a way that $a=(1\bis 1)$.
Now Malgrange's Preparation Theorem says
that there are smooth real-valued functions $v_\nu$ and $u_{\nu,j}$
defined in a neighborhood of $p=0$ such that
$$
\vp_\nu(x) = v_\nu(x) \left( x_1^k + 
\sum_{j=0}^{k-1}u_{\nu,j}(x')x_1^j \right)
$$
and $v_\nu(0)=1$.
The functions $u_{\nu,j}$ vanish of order $k-j$ at $x'=0$.
We see that in a neighborhood of $0$ the components $\vp_\nu$ have the
same zero set as the $x_1$-polynomials $x_1^k + 
\sum_{j=0}^{k-1}u_{\nu,j}(x')x_1^j$.
Lemma \ref{hyperfl} already tells us that $\NN(\vp_\nu)$
is countably $(n-1)$-$\Cu$-rectifiable near $0$.

If we denote the homogeneous part of degree $k-j$ of $u_{\nu,j}$
by $\hat{u}_{\nu,j}$, then we have
$$
\hat{\vp}_{p,\nu}(x) =  x_1^k + 
\sum_{j=0}^{k-1}\hat{u}_{\nu,j}(x')x_1^j  .
$$
By Lemma \ref{hauptteil} $\hat{\vp}_p$ solves
$$
\hat{L}_p \hat{\vp}_p = 0 .
$$
Since, by assumption, $\hat{\vp}_p\not= 0$ Lemma \ref{konstant} says
$\pi(\NN(\hat{\vp}_p)) \not= \R^{n-1}$.
Recall that $\pi:\R^n=T_pM \to \R^{n-1}$, $\pi(x)=x'$, is the projection
on the orthogonal complement of the $x_1$-axis.
Hence there exists an $x'_0$ such that the $x_1$-polynomials 
$x_1^k + \sum_{j=0}^{k-1}\hat{u}_{\nu,j}(x'_0)x_1^j$,
$\nu = 1 \bis N$, have no common root.
In other words, we can find linear combinations
\begin{eqnarray*}
\hat{F}(x) &=& \sum_{\nu=1}^N A_\nu \left( x_1^k + 
\sum_{j=0}^{k-1}\hat{u}_{\nu,j}(x')x_1^j \right) , \\ 
\hat{G}(x) &=& \sum_{\nu=1}^N B_\nu \left( x_1^k + 
\sum_{j=0}^{k-1}\hat{u}_{\nu,j}(x')x_1^j \right)
\end{eqnarray*}
which do not have a common root $x_1$ for $x'=x'_0$.

Recall that two polynomials have a common root if and only if their
resultant vanishes.
Let us denote the resultant of the two $x_1$-polynomials $\hat{F}$
and $\hat{G}$ by $R_{\hat{F},\hat{G}}(x')$.
It is a homogeneous polynomial in $x'$ with $R_{\hat{F},\hat{G}}(x'_0)\not=0$.
Moreover, $R_{\hat{F},\hat{G}}(x')$ is the lowest order term in the 
Taylor expansion of the resultant $R_{F,G}(x')$ of the $x_1$-polynomials
\begin{eqnarray*}
{F}(x) &=& \sum_{\nu=1}^N A_\nu \left( x_1^k + 
\sum_{j=0}^{k-1}{u}_{\nu,j}(x')x_1^j \right) , \\ 
{G}(x) &=& \sum_{\nu=1}^N B_\nu \left( x_1^k + 
\sum_{j=0}^{k-1}{u}_{\nu,j}(x')x_1^j \right) .
\end{eqnarray*}

If $x=(x_1,x') \in \NN(\vp)$ near $p$, then $F(x)=G(x)=0$ and hence
$R_{F,G}(x') = 0$.
In other words, near $p$, $\pi(\NN(\vp))\subset \NN(R_{F,G})$.
For each $x' \in \NN(R_{F,G})$ there are at most $k$ roots
$x_1$ such that $\vp(x_1,x')=0$.
>From the proof of Lemma \ref{hyperfl} we know that $\NN(\vp)$
is contained in a countable union of graphs of smooth functions
defined on the $x'$-hyperplane.
Hence $\NN(\vp)$ is countably $(n-2)$-$\Cu$-rectifiable if $\NN(R_{F,G})$ is.

>From $(\widehat{R_{F,G}})_0 = R_{\hat{F},\hat{G}} \not= 0$ we see that 
$R_{F,G}$ vanishes at $x'=0$ of finite order.
Applying Lemma \ref{hyperfl} once more we conclude that $\NN(R_{F,G})$
is countably $(n-2)$-$\Cu$-rectifiable.


\begin{center}
\pspicture(0,0)(14,7)

\psline(0,2)(4,2)
\psline(2,4)(6,4)
\psline(0,2)(2,4)
\psline(4,2)(6,4)

\rput(3,4){\psframebox*[framearc=0.2]{}}
\psline{->}(3,3)(3,6)
\psline[linestyle=dotted](3,3)(3,2)
\psline(3,2)(3,0.3)

\pscurve(1.5,3.5)(2.5,3.5)(3,3)(4,2.5)(4.5,2.5)
\pscurve(0.5,2.5)(2,2.5)(3,3)(4.5,3.5)(5.5,3.5)

\rput(3.5,5.5){$x_1$}
\rput(5.5,3.8){$x'$}
\rput(2,2.5){\psframebox*[framearc=0.5]{$\mathcal{N}(R_{F,G})$}}

\psline(6,2)(10,2)
\psline(8,4)(12,4)
\psline(6,2)(8,4)
\psline(10,2)(12,4)

\rput(9,4){\psframebox*[framearc=0.2]{}}
\psline{->}(9,3)(9,6)
\psline[linestyle=dotted](9,3)(9,2)
\psline(9,2)(9,0.3)

\pscurve(7.5,3.5)(8.5,3.5)(9,3)(10,2.5)(10.5,2.5)
\pscurve(6.5,2.5)(8,2.5)(9,3)(10.5,3.5)(11.5,3.5)

\rput(8.75,4){\psframebox*[framearc=0.2]{}}
\pscurve(7.5,5)(8.5,4.5)(9,3)
\pscurve[linestyle=dotted](9,3)(9.25,2.3)(9.5,2)
\psecurve(9.25,2.3)(9.5,2)(10,1.8)(10.5,1.7)(10.75,1.65)

\rput(9.93,4){\psframebox*[framearc=0.2]{}}
\pscurve(9,3)(10.5,4.5)(11.5,5)
\pscurve[linestyle=dotted](9,3)(8.5,2.5)(7.75,2)
\psecurve(8.5,2.5)(7.75,2)(7,1.62)(6.5,1.52)(6,1.47)

\rput(8.2,4.5){\psframebox*[framearc=0.5]{$\supset\mathcal{N}(\varphi)$}}
\rput(10.48,4){\psframebox*[framearc=0.5]{}}
\psline{->}(10.5,4.2)(10.5,3.7)
\rput(10.8,4){\psframebox*[framearc=0.5]{$\pi$}}

\endpspicture
\end{center}
\begin{center}
\bf Fig.~2
\end{center}

The resultant $R_{\hat{F},\hat{G}}$ is a homogeneous polynomial of degree
$k^2$ in $x_2 \bis x_n$.
By counting preimages and applying \cite[3.2.27]{federer69a} as in
the proof of Lemma \ref{hyperfl} we obtain the following bound on
the $(n-2)$-dimensional Hausdorff measure:
\begin{eqnarray*}
\HH^{n-2}(\NN(\vp) \cap W^{n}(0,r))  
&\le&
\eh n(n-1)\cdot k \cdot k^2 \cdot \vol(W^{n-2}(0,r)) \\
&=&
\eh n(n-1)k^3 (2r)^{n-2}.
\end{eqnarray*}
Hence
$$
\Theta^{\ast n-2}(\NN(\vp),p) \le
\frac{2^{n-3}n(n-1)k^3}{\alpha(n-2)}.
$$
\qed


\section{An Example with a Wild Zero Set}

In this section we want to give evidence that our Main Theorem is
more or less all that one can say about local properties of zero sets of 
solutions to first order elliptic equations.
The following theorem will show that the zero set can be extremely
irregular and that its Hausdorff dimension can be any {\em real}
number $d\in [0,n-2]$.

\begin{theorem}\punkt{}
\label{wild}
Let $n\ge 3$.
Let $A \subset \R^{n-2}$ be any closed subset.
Then there exists a linear elliptic differential operator $D$ defined
over $\R^n$ and a solution $\vp$ of $D\vp=0$ such that
$$
\NN(\vp) = \{ (0,0) \} \times A .
$$
Moreover, $D$ has the strong unique continuation property.
\end{theorem}

\noi
{\em Proof.}
Let $d$ be exterior differentiation of differential forms on $\R^n$
and let $\delta$ be its adjoint operator (with respect to the Euclidean
metric).
Then 
$$
D_1 = d+\delta : \Cu(\R^n,\bigoplus_{j=0}^n \Lambda^j T^\ast\R^n) \to
\Cu(\R^n,\bigoplus_{j=0}^n \Lambda^j T^\ast\R^n)
$$
is a Dirac operator.
Let us denote the coordinates of $\R^n$ by $(x,y,z)$ with $x,y \in \R$
and $z\in \R^{n-2}$.
Now $\vp_1 = xdy + ydx$ is a solution of $D_1\vp_1 = 0$ with
$\NN(\vp_1) = \{ (0,0) \} \times \R^{n-2}$.

Pick a smooth function $F:\R^{n-2} \to \R$ with $F^{-1}(0) = A$.
The map $\Psi : \R^n \to \R^n$, $\Psi(x,y,z) = (x,y-F(z),z)$, is a 
diffeomeorphism.
This diffeomorphism maps the linear subspace $\{ (0,0) \} \times \R^{n-2}$
to the submanifold $E=\{ (\xi,\eta,\zeta)\ |\ \xi=0, \eta = -F(\zeta) \}$.
By pulling back the Euclidean metric via $\Psi$ we obtain a new
Riemannian metric on $\R^n$.
Let $D_2$ be the corresponding Dirac operator with respect to this new metric,
i.e.\ $D_2 = d + \delta_\Psi$ where $\delta_\Psi$ is the adjoint 
operator of $d$ with respect to the new metric.
There exists an analogous solution $\vp_2$ of $D_2\vp_2 = 0$ with
$\NN(\vp_2) = E$.

Hence the operator
$$
D=
\left( 
\begin{array}{cc}
D_1 & 0 \\
0 & D_2
\end{array}
\right)
$$
acting on $\Cu(\R^n,\bigoplus_{j=0}^n \Lambda^j T^\ast\R^n \oplus
\bigoplus_{j=0}^n \Lambda^j T^\ast\R^n)$ has the solution
$\vp = (\vp_1,\vp_2)$ with zero set 
$$
\NN(\vp) = \NN(\vp_1) \cap \NN(\vp_2) = \{(0,0)\} \times A .
$$
Since both $D_1$ and $D_2$ are Dirac operators and hence have the
strong unique continuation property so has $D$.
\qed


\section{Scalar Elliptic Equations of Second Order}

In this section we study implications of our Main Theorem for scalar 
elliptic equations of second order.
For a smooth real-valued function $u:M\to\R$ the {\em critical zero set} 
is defined by
$$
\mathcal{N}_{crit}(u) = \{ x \in M\ |\ u(x)=0, du(x)=0 \}
$$
Obviously, $\NN(u) - \NN_{crit}(u)$ is a smooth hypersurface.

By a {\em scalar semilinear elliptic differential operator of second order}
we mean a map
$$
L : \Cu(M,\R) \to \Cu(M,\R)
$$
which is locally of the form
$$
Lu(x) = -\sum_{ij}a_{ij}(x)\frac{\partial^2u}{\partial x_i \partial x_j}(x)
+ F(x,u(x),du(x))
$$
where $(a_{ij})$ is a symmetric and positive definite matrix and $F$ is 
smooth and respects the zero section, $F(x,0,0)=0$.

\begin{theorem}\punkt
\label{scalar1}
Let $M$ be a connected $n$-dimensional 
differential manifold.
Let $L : \Cu(M,\R) \to \Cu(M,\R)$ be a scalar semilinear elliptic differential
operator of second order.
Let $u \in \Cu(M,\R)$, $u \not\equiv 0$, be a solution of $Lu=0$.

Then $\NN(u)$ is countably $(n-1)$-$\Cu$-rectifiable with locally finite
$(n-1)$-dimen\-sional Hausdorff measure.
Moreover, $\mathcal{N}_{crit}(u)$ is a countably $(n-2)$-$\Cu$-rectifiable set
with locally finite $(n-2)$-dimensional Hausdorff measure.
In particular, we have for the Hausdorff dimension
$$
\dim(\mathcal{N}_{crit}(u)) \le n-2
$$
and if $n=2$, then $\mathcal{N}_{crit}(u)$ is a discrete set.
\end{theorem}

This theorem has been shown by Hardt and Simon using different methods
in \cite[Thm. 1.10]{hardt-simon89a} except for the statement that
$\NN_{crit}$ has locally finite $(n-2)$-dimensional Hausdorff measure.
This latter statement has been shown for $n=3$ by M.\ and T.\ 
Hoffmann-Ostenhof and Nadirashvili \cite{hoffmann-ostenhof-nadirashvili96a}.
The case $n \ge 4$ has been settled very recently by Hardt \cite{hardt97ppa}.

\epnoi
{\em Proof.}
The principal symbol of $L$ defines a Riemannian metric on $M$ with
respect to which $L$ can be written as
$$
Lu = \Delta u + \tilde{F}(x,u,du)
$$
where $\Delta = \delta d$ is the Laplace-Beltrami operator and $\tilde{F}$
is smooth and respects the zero section.

Let $E = \bigoplus_{j=1}^n \Lambda^j T^\ast M$ be the exterior form bundle.
We define the smooth map $V : E \to E$ by
$$
V\left(\sum_{j=0}^n \omega_j\right) = 
\tilde{F}(x,\omega_0,\omega_1) - \omega_1
$$
where $\omega_j \in \Lambda^j T^\ast_x M$.
Note that $V$ respects the zero section.

Now if $u$ solves $Lu=0$, then $\omega = u + du \in \Cu(M,E)$ solves
$(d+\delta + V)(\omega) =0$.
The operator $D=d+\delta + V$ is a semilinear Dirac operator acting
on $\Cu(M,E)$.
Thus Corollary \ref{dirac} yields the statements on 
$\mathcal{N}_{crit}(u) = \mathcal{N}(\omega)$.

Since Dirac operators have the strong unique continuation property
Lemma \ref{hyperfl} shows that $\NN(u) = \NN_{fin}(u)$ has locally
finite $(n-1)$-dimensional Hausdorff measure.
\qed

\begin{theorem}\punkt
\label{scalar2}
Let $M$ be a closed $n$-dimensional Riemannian manifold.

Then there exists a constant $C = C(M)$ such that for any eigenfunction
$u \in \Cu(M,\R)$ of the Laplace operator, $\Delta u = \lambda u$, and 
any $p \in M$ we have
\begin{eqnarray*}
\Theta^{\ast n-1}(\NN(u),p) &\le& C \sqrt{\lambda} , \\
\Theta^{\ast n-2}(\NN_{crit}(u),p) &\le& C \lambda^{3/2} . 
\end{eqnarray*}
\end{theorem}

\noi
{\em Proof.}
A theorem of Donnelly and Fefferman \cite[Thm.~1.1]{donnelly-fefferman88a}
tells us that there is a constant $C'=C'(M)$ such that $u$ can vanish
only up to order $C'\sqrt{\lambda}$.
Now the Main Theorem and Lemma \ref{hyperfl} give the bounds on the
Hausdorff density.
\qed

\epnoi
This theorem gives us upper bounds on the Hausdorff measures of
$\NN(u)$ and $\NN_{crit}(u)$ in small balls.
Unfortunately, we have no explicit control on how small the radius of these
balls must be such that the estimates hold.
Therefore we cannot deduce bounds on $\HH^{n-1}(\NN(u))$ and 
$\HH^{n-2}(\NN_{crit}(u))$ in terms of a suitable power of $\lambda$.
It is an old conjecture of Yau \cite[Probl.~38]{yau93a} that
$$
\HH^{n-1}(\NN(u)) \le C(M) \sqrt{\lambda}.
$$
In the case of analytic Riemannian metrics this has been shown by Donnelly 
and Fefferman \cite[Thm.~1.2]{donnelly-fefferman88a} whereas in the smooth 
case it is still open.
See e.g.\ \cite[Thm.~5.3]{hardt-simon89a} and
\cite[Cor.~1.3]{donnelly-fefferman90a} for partial results.
See also \cite{han-hardt-lin97ppa} for a bound on the Hausdorff measure
of $\NN_{crit}(u)$.


\bibliographystyle{amsplain}
\bibliography{meine,allg}

\ep
\ep

Mathematisches Institut

Universit\"at Freiburg

Eckerstr.~1

79104 Freiburg

Germany

\ep

E-Mail:
{\tt baer@mathematik.uni-freiburg.de}

\end{document}